\documentclass[12pt]{article}
\usepackage{amssymb,amsmath}
\def\hang{\hangindent\parindent}
\def\tex#1{\indent\llap{[#1]\enspace}\ignorespaces}
\def\re{\par\hang\tex}

\def\a{\alpha}

\def\d{\delta}

\def\Ker{{\rm Ker}}

\def\si{\sigma}

\def\ssc{\scriptscriptstyle}

\def\cl{\centerline}

\def\rar{\rightarrow}

\def\vs{\vspace*}

\def\ni{\noindent}
\def\ptl{\partial}

\def\Z{\mathbb{Z}{\ssc\,}}

\def\C{\mathbb{C}{\ssc\,}}
\def\F{\mathbb{F}{\ssc\,}}

\def\QED{\hfill$\Box$} \numberwithin{equation}{section}
\newtheorem{theo}{Theorem}[section]
\newtheorem{defi}[theo]{Definition}
\newtheorem{rema}[theo]{Remark}
\newtheorem{lemm}[theo]{Lemma}
\newtheorem{coro}[theo]{Corollary}
\newtheorem{prop}[theo]{Proposition}

\def\Id{{\rm Id}}
\begin{document}

\cl{{\Large \bf
 Simple deformed Witt algebras }\footnote{Supported by NNSF grant 10871125, 11071147 of China, NSF grant Y2006A17, Y2008A04 of Shandong Province,
   NSF J06P52 of Shandong Provincial Education Department, China }}
\vs{9pt}

\cl{ Guang'ai Song$^*$, Chunguang Xia$^\dag$}

\cl{\small $^{*\,}$College of Mathematics and Information Science,
Shandong Institute of \vs{-4pt}Business }

\cl{\small and Technology, Yantai, Shandong 264005, China}

\cl{\small $^{\dag\,}$Department of Mathematics, University of
Science and Technology of \vs{-4pt}China}\cl{\small Hefei 230026,
China }

\cl{\small E-mail: gasong@ccec.edu.cn, chgxia@mail.ustc.edu.cn}
\vs{15pt}

{\baselineskip0pt\lineskip0pt\begin{minipage}{5.8truein}
\baselineskip0pt\lineskip0pt\noindent{\small{\bf Abstract.} For any
factorization domain $\cal A$ and an algebra endomorphism $\sigma$
of $\cal A$, there exists a non-associative algebra $({\cal
A},\sigma,[\cdot,\cdot])$ with multiplication satisfying
skew-symmetry and generalized (twisted) Jacobi identities, called a
$\sigma$-deformed Witt algebra. In this paper, we obtain the
necessary and sufficient conditions for the algebra $({\cal
A},\sigma,[\cdot,\cdot])$ to be simple. \vskip10pt

  }

 \noindent{\bf 2000 MR Subject Classification:} 17B66, 17B37
\end{minipage}}\vskip20pt

\cl{\bf\S1. \
Introduction}\setcounter{section}{1}\setcounter{equation}{0}\vskip7pt

Deformation theory has many applications in mathematics and physics,
especially in quantum theory (e.g., [D, ES, Kas]), therefore more
and more attentions have been received on this field. For instance,
the $q$-Witt algebras, $q$-Virasoro algebra, $q$-deformed Heisenberg
algebras were investigated in [H, HS], and the quantized structures
of some Witt type algebras were present in [G, HW]. Furthermore,
there are various interesting deformations of Lie algebras, such as
hom-Lie algebras, quasi-Lie algebras, quasi-hom-Lie algebras (e.g.,
[LS1, LS2, LS3, MS, RS, Y]). As a by-product, some algebras such as
hom-Lie color algebras, hom-admissible algebras, hom-coalgebras,
hom-Hopf algebras, etc. have also received some attentions in
literature  (e.g., [HLS, HS, LS1, LS2, LS3, MS, Y]).

The Witt algebra, which is the infinite-dimensional Lie algebra of
linear differential operators on the circle, is as stated in [HLS]
an important example in the classical differential and integral
calculus, relating it to topology and geometry, and at the same time
responsible for many of its key algebraic properties. It is also a
well-known fact that the Virasoro algebra, which has played more and
more important role in various mathematical and physical theories,
such as vertex operator algebras, conformal field theory, the theory
of the quantum Hall effect and integrable system, etc, is simply the
universal central extension of the Witt algebra. Furthermore, the
Witt type Lie coalgebras are the first examples of infinite
dimensional Lie coalgebras constructed by Michaelis, and called
``the gift given by God'' [M]. Due to the increasingly importance of
Witt type algebras, generalizations of Witt type algebras and their
subalgebras have been introduced and studied by many authors (e.g.,
[DZ1, DZ2, Kaw, O, P, S, SX1, SX2, SX3, SXZ, SZ, SZh, X]). In
particular, the necessary and sufficient conditions for a wide class
of Witt type algebras to be simple were obtained in [P], which can
be conveniently used to verify whether or not a Witt type algebra is
simple.
 From this, a large class of simple
Witt type algebras were explicitly constructed in [X], whose
structure  and representation theories were studied in [SXZ, SZ,
SZh].

Using $\si$-derivations, the authors in [HLS] obtained general
deformations of Witt type algebras, referred to as {\it
$\si$-deformed Witt algebras.} In order to construct explicitly
simple $\si$-deformed Witt algebras and to study their structure and
representation theories, it is desirable to establish some kind of
criteria similar to [P] for Witt type algebras, to determine when a
$\si$-deformed Witt algebra is simple. This is the aim of the
present paper.

The main result of the present paper is the following theorem.
\begin{theo}\label{theo-1}
 Let $\cal{A}$ be a unique factorization domain over an arbitrary field $\F$,
 $\si$ an epimorphism on $\cal{A}$, and $\ptl$ the
$\si$-derivation defined by $(\ref{to-def})$ satisfying
$\ptl(\cal{A}) = \cal{A}$. Then the $\si$-deformed Witt  algebra
$(\cal{A}, [\cdot, \cdot], \si )$ defined in Theorem $\ref{theo-2}$
is simple if and only if $\cal{A}$ is $ \ptl$-simple
$($cf.~~Definition $\ref{stable})$.
\end{theo}
\vskip10pt

\cl{{\bf\S2.  Preliminary
results}}\vskip7pt\setcounter{section}{2}\setcounter{theo}{0}\setcounter{equation}{0}
 Let $\F$ be an arbitrary field, $\cal{A}$ a unique factorization domain over $\F$,
 and $\si$ an endomorphism on $\cal{A}$. we first recall the notion of $\si$-derivation of $\cal{A}$(see [Kas]),
 A linear map $d: \cal{A}
\rar \cal{A}$ is called a {\it $\si$-derivation of $\cal{A}$} if
\begin{equation}\label{si-derv}
d(ab) = d(a)b + \si(a) d(b)\mbox{ \  for all }a, b \in \cal{A}.
\end{equation}
Denote the vector space of $\si$-derivations of $\cal{A}$ by
$\cal{D}(A)$. It is easy to see that $d(1) = 0$ for $d \in
\cal{D}(A)$.

The following notions can be found in [HLS].  Assume that
$\si\ne\Id$ (the identity map). Let $g = \mbox{\rm gcd}(\Id -
\si)({\cal A})$ be a greatest common divisor of $\,(\Id - \si)({\cal
A})$.  Define a linear map $\ptl:{\cal A}\to{\cal A}$ by
\begin{equation}
\ptl = \frac{\Id - \si}{g }, \mbox{ \ i.e., } \ptl(x) = \frac{x -
\si (x)}{g} \mbox{ \ for }x \in \cal{A}.\label{to-def}
\end{equation}
The following result can be found in [HLS].

\begin{lemm}\label{lemm-1}
 The space
$\cal{D}(\cal{A})$ is free of rank one as an $\cal{A}$-module with
generator $\ptl$.
\end{lemm}

 Thus ${\cal D}(\cal{A}) =
{\cal{A}}\ptl$, where an element  $a \ptl \in \cal{A} \ptl$ acts on
$\cal{A}$ as $$a \ptl: x \rar a \ptl(x) \mbox{ \ for  }a, x \in
\cal{A}.$$ The algebra $(\cal{A}, \si, [\cdot,\cdot])$ with product
defined below, is called a \textit{$\si$-deformed Witt algebra}:
\begin{equation}
[a, b] = \si(a)\ptl(b) - \si(b)\ptl(a)\mbox{ \ for \ }a,b\in{\cal
A}.\label{the-2}
\end{equation}
The following  is one of the main results in [HLS].

\begin{theo}\label{theo-2}
\noindent The $\si$-deformed Witt algebra $(\cal{A}, \si,
[\cdot,\cdot])$
 satisfies the following conditions:
\begin{itemize}\parskip-3pt\item[{\rm(1)}] bilinearity: The operator $[\cdot, \cdot]: \cal{A}
\times \cal{A} \rar \cal{A}$ defined by $(\ref{the-2})$ is bilinear.

\item[{\rm(2)}] skew-symmetry: $[a, a ] =0 {\rm \ for \ all} \ \ a\in \cal{A}.$

\item[{\rm(3)}] the generalized Jacobi identity:
$$\begin{array}{ccl}0 &=&[\si(a), [b, c ]] + \d [a , [b, c ]]\\
&&+ [\si(b), [c, a ]] + \d [b, [c, a ]]\\
&&+ [\si(c), [a, b ]] + \d [c, [a, b ]]\mbox{ \ for }a, b, c \in
\cal{A},
\end{array}$$
where, $\d = \frac{\si(g)}{g}.$\end{itemize}
 \end{theo}
\ni {\it Proof.~}~The proof is straightforward.\hfill$\Box$

\begin{defi}\label{hom} A Hom-Lie algebra is a triple $(V, [\cdot, \cdot],
\a)$ consisting of a vector space $V$, bilinear map $[\cdot, \cdot]:
V \times V \rar V$ and a linear map $\a: V \rar V$
satisfying \begin{eqnarray*}&& [x, y] = -[y, x], \\[4pt]
&& [\a(x), [y, z]] + [\a(y), [z, x]] + [\a(z), [x, y]]
=0,\end{eqnarray*} for all $x, y, z \in V.$
\end{defi}
 \begin{rema}\label{rema1} \rm \begin{itemize}\parskip-3pt\item[(1)]
 Since $g| (\Id -\si) (g)$ from the definition of $g$,
we have $g | \si(g),$ and so $\d = \frac{\si(g)}{g} \in \cal{A}.$
Thus the generalized Jacobi identity is well defined in $\cal{A}$.
Furthermore, $\d$ satisfies: $$\ptl(\si(a)) = \frac{(\Id -
\si)(\si(a))}{g} = \frac{\si(a)- \si^2(a)}{g} = \frac{\si(g)}{g}
\si(\frac{a - \si(a)}{g}) = \d \si(\ptl(a)),$$ for  $a \in \cal{A}.$
\item[(2)] If $\si \rar 1$ (in some sense),  the $\si$-derivation $\ptl$ defined by
(\ref{to-def}) is the usual derivation on $\cal{A}$, and the bracket
defined by (\ref{the-2}) is the bracket of usual Witt type Lie
algebras.
 \item[(3)] 
If $\d\in\F$,
then by defining a linear map $\si_1: \cal{A} \rar \cal{A}$ by
$\si_1(a) = \si(a) + \d a$ for $a \in \cal{A}$, and the bracket on
$\cal{A}$ by $[a, b] = \si(a)\ptl(b) - \si(b)\ptl(a)$ for $a, b \in
\cal{A},$ we obtain a Hom-Lie algebra $({\cal{A}}, [\cdot, \cdot],
\si_1)$, called a {\it Witt-type Hom-Lie algebra}.
\item[(4)] If $\d \notin \F$,  the triple $({\cal{A}}, [\cdot, \cdot],
\si_1)$ is not a Hom-Lie algebra, in this case, the the triple
$({\cal{A}}, [\cdot, \cdot], \si_1)$ is no longer satisfying
Hom-Jacobi identity defined by Definition \ref{hom}, and one can
compute it directly.
\end{itemize}\end{rema}\vskip10pt

\cl{{\bf\S3. The Proof of Main Result
}}\vskip7pt\setcounter{section}{3}\setcounter{theo}{0}\setcounter{equation}{0}

From now on, we shall discuss the conditions for $\si$-deformed Witt
algebra $({\cal{A}}, [\cdot, \cdot], \si)$ to be simple. We shall
always assume $\cal A$ is a unique factorization domain. As in the
Lie algebra case (e.g., [P]), we first give some definitions.

\begin{defi}\rm\label{stable} \begin{itemize}\parskip-3pt\item[(1)]An ideal $I$ of $\cal{A}$ is {\it
$(\si, \cal{D}(\cal{A}))$-stable} if $ \si(a), d(a) \in I$ for all
$a \in I,\, d \in {\cal D}(\cal{A}).$ \item[(2)]The algebra
$\cal{A}$ is {\it $(\si, \cal{D(\cal{A})})$-simple} if it has no
nontrivial $(\si, \cal{D(\cal{A})})$-stable ideals.
\end{itemize}\end{defi}
\begin{rema}\label{rema-2} \rm\begin{itemize}\parskip-3pt\item[(1)]It is easy to see that if $I$ is a $(\si,
\cal{D(\cal{A})})$-stable ideal of $\cal{A}$, then $I$ is also an
ideal of the $\si$-deformed Witt algebra $(\cal{A}, [\cdot, \cdot],
\si).$
\item[(2)]
By Lemma \ref{lemm-1}, $\cal{D}(\cal{A}) = \cal{A}\ptl$,  thus an
ideal $I$ of $\cal{A}$ is $\cal{D}(\cal{A})$-stable if and only if
it is $\ptl$-stable.
\end{itemize}\end{rema}

Similar to the usual Witt-type Lie algebra case (see [P]), we have

\begin{prop} If $\cal{A}$ is $(\si,
\cal{D(\cal{A})})$-simple, then $\cal{A}^{\ptl}$ $= \{ a\, |\,
\ptl(a) = 0, a \in $ $\cal{A} \}$ is $(\si, \ptl)$-stable, moreover,
$\cal{A}^{\ptl}$ is a field containing $\F.$
\end{prop}
\ni{\it Proof.~}~From the definition of $\ptl,$ we see that $\ptl(x)
= 0$ if and only if $\si (x) =x.$ Thus $\cal{A}^{\ptl}$ is obviously
$\si$-stable. It is easy to see that $\cal{A}^{\ptl}$ is a subring
of $\cal{A}$. If $0 \neq a \in \cal{A}^{\ptl},$ then $a \cal{A}$ is
a nonzero $(\si, \cal{D(\cal{A})})$-stable ideal of $\cal{A},$ and
so $a \cal{A} = \cal{A}.$ Thus $a$ is an invertible element in
$\cal{A}.$ From $0= \ptl(1) = \ptl(a a^{-1}) = \ptl(a) a^{-1} +
\si(a) \ptl(a^{-1}) = a \ptl(a^{-1}),$ and since $\cal{A}$ has no
nonzero divisors, we have $\ptl(a^{-1}) = 0.$ Hence $a^{-1} \in
\cal{A}^{\ptl},$ namely, $\cal{A}^{\ptl}$ is a field.
 \QED
\begin{lemm}\label{lemm-2} 
If an ideal $I$ of $\cal{A}$ is $\ptl$-stable, then $I$ is also
$\si$-stable.
\end{lemm}
\ni{\it Proof.~}~For any $x \in I, $ since $ \frac{x -
\si(x)}{g}=\ptl(x)\in I$, we obtain $\si(x) = x - g\ptl(x) \in I.$
\QED\vskip7pt

From Lemma \ref{lemm-2} and Remark \ref{rema-2}, we have,
\begin{eqnarray}\label{aS1}
\mbox{ an ideal $I$ of $\cal{A}$ is $(\si, {\cal
D}(\cal{A}))$-stable}&
\Longleftrightarrow&\mbox{$I$ is $\ptl$-stable,}\\
\label{aS2}\mbox{$\cal{A}$ is $(\si, {\cal
D}(\cal{A}))$-simple}&\Longleftrightarrow&\mbox{$\cal{A}$ is
$\ptl$-simple.}
\end{eqnarray}
\begin{lemm}
If $\si \neq 0$ and $\cal{A}$ is $ \ptl$-simple, then $\si$ is a
monomorphism.
\end{lemm}
\ni{\it Proof.~}~Suppose that $\si$ is not a monomorphism, then
$\Ker(\si) \neq \{0\}$ is an ideal of $\cal{A},$ which is obviously
$\si$-stable. For any $x \in \Ker(\si),$ we have $\ptl(x) =
\frac{(\Id - \si)(x)}{g} = \frac{x}{g},$ and $\si(\ptl(x)) =
\si(\frac{x}{g}) = \frac{\si(x)}{\si(g)} = 0,$ and so $\ptl(x) \in
\Ker(\si)$. Thus $\Ker(\si)$ is a nontrivial $(\si, \ptl)$-stable
ideal of $\cal{A}$ by (\ref{aS1}), a contradiction. \QED

\begin{lemm}\label{lemm-3} Let $V$ be a $\ptl$-stable vector subspace of
$\cal{A}$, and $\si$ an epimorphism of $\cal{A}$. Then any maximal
ideal $I$ of $\cal{A}$ which is contained in $V$ is $(\si, {\cal
D}(\cal{A}))$-stable.
\end{lemm}
\ni{\it Proof.~}~It is easy to see that $I + \ptl(I) \subset V$ is
$\ptl$-stable. For $x \in I,\, a \in \cal{A}$, from $\ptl(ax) =
\ptl(a) x + \si(a) \ptl(x)$, we have $\si(a)\ptl(x) = \ptl(ax) -
\ptl(a) x \in I + \ptl(I)$. Since $\si$ is an epimorphism of
$\cal{A},$ we see that $I +\ptl(I)$ is an ideal of $\cal{A}$. So $I
+ \ptl(I) \subset I$ by the maximal nature of $I.$ By (\ref{aS1}),
$I$ is $(\si, {\cal D}(\cal{A}))$-stable. \QED

\begin{lemm} \label{lemm-3.8} Suppose  $\cal{A}$ is $ \ptl$-simple. Let $I$
be an ideal of $(\cal{A}, [\cdot, \cdot], \si)$, and $\si$ an
epimorphism on $\cal{A}.$ Then $\ptl(\cal{A}) $ $\subset I.$
\end{lemm}
\ni{\it Proof.~}~For $x\in I$, we have $\ptl(x)=
\si(1)\ptl(x)-\si(x)\ptl(1)=[1,x]\in I,$ i.e., $I$ is a
$\ptl$-stable subspace of $\cal{A}$.
For  $x \in I, a \in \cal{A},$ we have
$
\si(x)\ptl(a) -
\si(a)\ptl(x)=[x,a] \in I . $ 
\noindent From this, we
obtain
\begin{equation*}
\ptl(xa) =\ptl(x)a + \si(x)\ptl(a) \equiv \ptl(x) a + \si(a)\ptl(x)
\ ({\rm mod\,} I).
\end{equation*}
Thus $\ptl(xa) \in I\cal{A},$ i.e., $I\cal{A}$ is a $\ptl$-stable
ideal of $\cal{A}.$

Let $x$ be an element of $I$ with the shortest expression: $x = a_x
x_1 x_2 \cdots x_s,$ where $a_x$ is a unit of $\cal{A}$, and $x_i,
i=1, 2, \cdots, s$ are irreducible elements of $\cal{A}$ with $s$
being minimal. Then $x\cal{A}$ must be a maximal ideal of $\cal{A}$
which is contained in $I\cal{A}$. Thus $x \cal{A}$ is a
$\ptl$-stable ideal of $\cal{A}$ by Lemma \ref{lemm-3}. This implies
that $x\cal{A} = \cal{A}$  since $\cal{A}$ is $\ptl$-simple.
Therefore $x$ is a unit of $\cal{A}.$ For any $a\in{\cal A}$, we
have $$\begin{array}{ll}\si(x)x\ptl(a)\!\!&=\si(x)(\ptl(a)x +
\si(a)\ptl(x)) - \si(a)\si(x)\ptl(x)\\[4pt]&=\si(x)\ptl(ax) -
\si(ax)\ptl(x)\\[4pt]&=[x,ax]\in I.\end{array}$$ This implies $\ptl(a) \in I$
 since $\si(x)x$
is also a unit. This proves $\ptl(\cal{A})$ $ \subset I.$
 \QED\vskip7pt

Now we can prove the main result of the paper.\vskip7pt

\ni {\it Proof of Theorem $\ref{theo-1}$.~}~On one hand, if
$\cal{A}$ is not $ \ptl$-simple, and $I$ is a $\ptl$-stable proper
ideal of $\cal{A},$ then one can easily check that $(I, [\cdot,
\cdot], \si)$ is a proper ideal of $(\cal{A}, [\cdot, \cdot], \si).$

 On the other hand, if $(I, [\cdot, \cdot], \si)$ is an ideal of $(\cal{A}, [\cdot, \cdot], \si)$,
 then
${\cal A}=\ptl({\cal A})\subset I$ by Lemma \ref{lemm-3.8}.
\hfill$\Box$\vskip7pt

From Theorems \ref{theo-1} and  \ref{theo-2}, we obtain
 \begin{theo}\label{theo-3} Let $\cal{A}$ be a unique factorization
domain over $\F$, $\si$  an epimorphism on $\cal{A}$, and $\ptl$ the
$\si$-derivation defined by {\rm (\ref{to-def})} satisfying
$\ptl(\cal{A}) = \cal{A}$. Suppose  $\d$ defined in Theorem {\rm
\ref{theo-2}} is an element in $\F$. Then the Hom-Lie algebra
$({\cal A}, [\cdot, \cdot],\si_1)$ with $\si_1$ being defined in
Remark {\rm \ref{rema1}(3)} is simple if and only if $\cal{A}$ is $
\ptl$-simple.
\end{theo}\vskip10pt

\cl{{\bf\S4. Applications
}}\vskip7pt\setcounter{section}{4}\setcounter{theo}{0}\setcounter{equation}{0}
As applications of Theorems \ref{theo-1} and \ref{theo-3}, we obtain
the following.

 \begin{coro}\label{prop-4.1} Let ${\cal A} = \C[t]$ be the polynomial
algebra in one variable t, and $\si$ the endomorphism of $\cal A$
determined by $ \si (t) = qt,$ where $0\ne q\in\C$ is not a root of
unit. Thus $$\si(f(t)) = f(qt)\mbox{ \ for \ }f(t) \in \cal{A}.$$
Define the $\si$-derivation by \begin{equation}\label{4---1}\ptl
(f(t)) = \frac {(\Id - \si)(f(t))}{t-qt} = \frac{f(t) -
f(qt)}{t-qt}.\end{equation}Then the $\si$-deformed Witt algebra
$(\cal{A}, [\cdot, \cdot], \si)$ is simple. Furthermore, the Hom-Lie
algebra $({\cal{A}}, [\cdot, \cdot], \si_1)$ with  $\si_1$ being
defined in Remark {\rm \ref{rema1}(3)}  is simple.
\end{coro}
\ni{\it Proof.~}~It easy to check that $g = \gcd (\Id -
\si)({\cal{A}}) = t - qt,$ and $\ptl(\cal{A}) = \cal{A}.$ Let $I
\neq \{0\}$ be a $\ptl$-stable ideal of $\cal{A}.$ Since $\cal{A}$
is a principal ideal domain, we have $I = (p(t))$ for some $0 \neq
p(t) \in \cal{A}.$ Note that $\ptl(p(t))\in I=(p(t))$, i.e.,
$p(t)|\ptl(p(t))$. On the other hand, we have
$\deg\,\ptl(p(t))<\deg\,p(t)$ by (\ref{4---1}). This implies
$\ptl(p(t))=0$, namely, $p(t)\in\C$ is invertible. Thus $I =
\cal{A},$ and
 $\cal{A}$ is $\ptl$-simple. By Theorem \ref{theo-1},
 $(\cal{A}, [\cdot, \cdot], \si)$ is simple.

  Since $\d =
\frac{\si(g)}{g} = \frac{\si(t) - \si(qt)}{t-qt} =q \in \F, $ by
Theorem \ref{theo-3}, $(\cal{A}, [\cdot, \cdot],
$$\si_1)$ is a simple Hom-Lie algebra, where $\si_1$ is defined as in
Remark \ref{rema1}(3) by $\si_1 = \si + q\, {\rm \Id}.$\QED

\begin{coro} Let $\cal{A}$, $\si, {\rm and} \ \ptl$ be as in Corollary {\rm \ref{prop-4.1}}, and
$q$ be an $n$-th primitive root of unit. Then the $\si$-deformed
Witt algebra $(\cal{A}, [\cdot, \cdot], \si)$ is not simple.
\end{coro}
\ni{\it Proof.~}~The principal ideal $I=(t^n)$  generated by $t^n$
is a proper $(\si, \ptl)$-stable ideal of $\cal{A}$, and  $(I,
[\cdot, \cdot], \si)$ is a proper ideal of $(\cal{A}, [\cdot,
\cdot], \si)$. Thus $(\cal{A}, [\cdot, \cdot], \si)$ is not simple.
\QED

\begin{coro} Let ${\cal A} = \C[t^{\pm 1}]$ be the Laurent polynomial
algebra in one variable $t,$  and $\si$ the endomorphism on
$\cal{A}$ determined by $\si(t) =qt$, where $0\ne q\in\C$ is not a
root of unit. Then the Hom-Lie algebra $({\cal A}, [\cdot, \cdot],
\si_{1})$ with $\si_1$ being defined in  Remark {\rm \ref{rema1}(3)}
is simple.
\end{coro}
\ni{\it Proof.~}~Since  $(\Id - \si)(t) = t - qt = t (1-q)$ is a
unit in $\cal{A}$, we see that $g=\gcd(\Id-\si)({\cal A})$ is a unit
of $\cal A$.

 Take $g=t^k$ for any fixed $k\in\Z$. Then $\ptl = t^{-k}(\Id - \si)$.
 It is easy to see
that $\ptl(\cal{A}) = \cal{A}.$ Let $\cal{I}$ be a nonzero
$\ptl$-stable ideal of $\cal{A}$ (thus also $\si$-stable). Let $0\ne
p(t)=\sum_{i=-m}^na_it^i\in{\cal I}$. By applying $\si$ to $p(t)$
several times, we obtain
\begin{equation}\label{Asystem}\sum_{i=-m}^nq^{ij}(a_it^i)=\si^j(p(t))\in{\cal
I},\ \ \ j=0,1,2,...,m+n.\end{equation} Regarding (\ref{Asystem}) as
a system of linear equations on $(m+n+1)$ variables
$a_it^i,\,i=-m,-m+1,...,n$, since the determinant of coefficients is
a nonzero Vandermonde determinant due to the fact that $q$ is not a
root of unit, we obtain that each monomial (which is a unit in $\cal
A$) of $p(t)$ is also in $\cal I$. Thus $\cal{I } = \cal{A}.$ Let
$\d = \frac{\si(g)}{g} = q^{k},$ and $\si_1 = \si + q^k \,\Id,$ we
obtain from Theorem \ref{theo-3} that $(\cal{A}, [\cdot, \cdot],
$$\si_1)$ is a simple Hom-Lie algebra. \QED

\begin{coro} Let ${\cal A} = \C[t^{\pm 1}]$ be the Laurent polynomial
algebra in one variable $t,$ and $\si$ the endomorphism of $\cal{A}$
defined by $\si(t) = q t^s$, where $0\ne q\in\C,\,s \in \Z,\, s \neq
0, 1,2$. Then the $\si$-deformed Witt algebra $({\cal A}, [\cdot,
\cdot], \si)$ defined by
 Theorem {\rm \ref{theo-2}} is not simple. Moreover, the triple $({\cal A}, [\cdot, \cdot],
 \si_1)$ is no longer a Hom-Lie algebra.
\end{coro}
\ni{\it Proof.~}~It is easy to see that $g$ is determined by the
images of $t, t^{-1}$ on $\Id - \si.$ Since $$(\Id - \si)(t) = t -
qt^s = t (1 - q t^{s-1}), \ \ (\Id - \si)(t^{-1}) = t^{-1} -
q^{-1}t^{-s} = -q^{-1}t^{-s}( 1- q t^{s -1}),$$  and $ t,
-q^{-1}t^{-s}$ are units of $\cal{A}$, we can take $g =1 - q
t^{s-1}$ if $s>2$, or $g=1 - q^{-1} t^{1-s}$ if $s<0$.

 If $s > 2,$ by denoting $T = q t^{s-1},$ we have $$\ptl(T) = \frac{T- \si(T)}{1 - T}
 = \frac{T - T^s}{1- T} = T(1 + T + T^2 + \cdots + T^{s-2} ).$$
 Hence
 the ideal $I$ generalized by $1 + T + \cdots + T^{s-2}$ is a proper ideal
  of $\cal{A}$, which is easily checked to be $\ptl$-stable (thus also $\si$-stable by Lemma
 \ref{lemm-2}). Therefore $I$ is a proper ideal of
  $\si$-deformed Witt algebra $(\cal{A}, [\cdot, \cdot], \si)$, and so $(\cal{A}, [\cdot, \cdot],
  \si)$ is not simple.

  If $s<0,$ by denoting $T = q^{-1} t^{1-s}$, we have
  $$\ptl(T) = \frac{T- \si(T)}{1-T} = \frac{T-T^s}{1 - T}=
  -T^s\frac{1 - T^{1-s}}{1 -T} = -T^s (1 + T + \cdots + T^{-s}).$$
As above, the proper ideal $I$ of $\cal{A}$ generalized by $1 + T +
\cdots + T^{-s}$ is $(\si, \ptl)$-stable, so $I$ is an ideal of
$(\cal{A}, [\cdot, \cdot], \si),$ and the algebra $(\cal{A}, [\cdot,
\cdot], \si)$ is not simple.

In both cases, since $\frac{\si(g)}{g} \notin \F,$ the triple
$({\cal{A}}, [\cdot, \cdot], \si_1)$ is not a Hom-Lie algebra by
Remark \ref{rema1}(4). \QED

\begin{coro} Let ${\cal A}=\C[x_1^{\pm1}, x_2^{\pm1}, \cdots, x_n^{\pm1}]$ be the Laurent
polynomial with $ n$ variables, and $\si$ the endomorphism of
$\cal{A}$ defined by $\si(x_i) = q_i x_i,\, i = 1, 2, \cdots, n,$
where every $0\ne q_i\in\C$ is not a root of unit. Then the triple
$({\cal{A}}, [\cdot, \cdot], \si)$ is a simple $\si$-deformed Witt
algebra, and the triple $({\cal{A}}, [\cdot, \cdot], \si_{1})$ is a
simple Hom-Lie algebra.
\end{coro}
\ni {\it Proof.~}~It easy to check that $\si$ is an epimorphism of
$\cal{A}$, and we can take $g =\gcd(\Id-\si)({\cal A})=1.$ For any
monomial $x_1^{k_1} x_2^{k_2} \cdots x_n^{k_n} \in \cal{A}$, since
$\ptl(x_1^{k_1} x_2^{k_2} \cdots x_n^{k_n}) = (1 -
q_1^{k_1}q_2^{k_2} \cdots q_n^{k_n}) x_1^{k_1} x_2^{k_2} \cdots
x_n^{k_n},$ we see $\ptl(\cal{A}) = \cal{A}.$ Let $I$ be any nonzero
$\ptl$-stable ideal of $\cal{A}.$ As in (\ref{Asystem}), we can
prove $I={\cal A}$. Thus $\cal{A}$ is $\ptl$-simple, and the triple
$(\cal{A}, [\cdot, \cdot], \si)$ is a simple $\si$-deformed Witt
algebra. Since $\d =1$, by Theorem \ref{theo-3}, we obtain that the
triple $(\cal{A}, [\cdot, \cdot], \si_{\rm 1})$ is a simple Hom-Lie
algebra with $\si_1 = \si + \Id $. \QED\vskip15pt

\cl{\bf References}\vskip7pt \small
\parindent=8ex\parskip=6pt\baselineskip=3pt

\re{D} V. Drinfel'd, {Quantum groups}, in ``Proceeding of the
International Congress of Mathematicians",  Vol. 1, 2, Berkeley,
Calif. 1986, Amer. Math. Soc., Providence, RI, 1987, pp. 789--820.

\re{DZ1} D.$\check{\rm Z}$. Dokovi\'c,
  K. Zhao, {Generalized Cartan type W Lie algebras in characteristic
0,} {\it J. Algebra} {\bf195} (1997), 170--210.

\re{DZ2} D.$\check{\rm Z}$. Dokovi\'c,
  K. Zhao, isomorphisms,
and second cohomology of generalized Witt algebras, {\it Trans.
Amer. Math. Soc.} {\bf350} (1998), 643--664.

\re{ES} P. Etingof, O. Schiffmann, {\it Lectures on Quantum Groups},
2nd, International Press, USA, 2002.

\re{G} C. Grunspan, {Quantizations of the Witt algebra and of simple
Lie algebras in characteristic $p$}, {\it J. Algebra} {\bf280}
(2004), 145--161.

 \re{H} N. Hu, {$q$-Witt Algebras, $q$-Virasoro Algebra, $q$-Lie
Algebras, $q$-Holomorph Structure and Representations,} {\it Algebra
Colloq.} {\bf 6} (1999),No. 1, 51-70.

\re{HLS}J.T. Hartwig,  D. Larsson, S.D. Silvestrov, {Deformations of
Lie algebras using $\si$-derivations,} {\it J. Algebra} {\bf 295}
(2006), no. 2, 314--361.

\re{HS} L. Hellstr$\ddot{o}$m, S.D. Silvestrov,  {\it Commuting
elements in $q$-deformed Heisenberg algebras}, World Scientific
Publishing Co., Inc., River Edge, NJ, 2000.

\re{HW} N. Hu, X. Wang, {Quantizations of generalized-Witt algebra
and of Jacobson-Witt algebra in modular case}, {\it J. Algebra}
{\bf312} (2007), 902--929.


 \re{Kas} C. Kassel, {\it Quantum Groups}, Graduate Texts in
Mathematics, {\bf 155}, Springer-Verlag, New York, 1995.

\re{Kaw} N. Kawamoto,{ Generalizations of Witt algebras over a field
of characteristic zero,} {\it Hiroshima Math. J.} {\bf16} (1986),
417--462.

\re{LS1} D. Larsson, S.D. Silvestrov, {Qusi-hom-Lie algebras,
Central extensions and 2-cocycle-like identities,} {\it J. Algebra}
{\bf 288} (2005), 321--361.

\re{LS2}D. Larsson, S.D. Silvestrov, {Quasi-Lie algebras,} in
``Noncommutative Geometry and Representation Theory in Mathematical
Physics", {\it Contemp, Math.} {\bf 391} (2005), Amer. Math. Soc.,
Providence, RI, 241--248.

\re{LS3}D. Larsson, S.D. Silvestrov, {Quasi-deformations of $sl_2
(\F)$ using twisting derivations,} {\it Comm. Algebra} {\bf 35}
(2007), 4303--4318.

\re{M} W. Michaelis, {A Class of Infinite-dimensional Lie Bialgebras
Containing the Virasoro Algebras, {\it Adv. Math.} {\bf107} (1994),
365--392.

\re{MS}A. Makhlouf,   S.D. Silvestrov, {Hom-algebra structures,}
{\it J. Gen. Lie Theory Appl.} {\bf2} (2008), 51--64.
%

\re{O} J.M. Osborn, {New simple infinite-dimensional Lie algebras of
characteristic 0,} {\it J. Algebra} {\bf185} (1996), 820--835.

\re{P} D.P. Passman, Simple Lie algebras of Witt Type, {\it J.
Algebra} {\bf206} (1998), 682--692.

\re{RS} L. Richard, S.D. Silvestrov, {Qusi-Lie structure of
$\si$-derivations of $\C[t^{\pm1}]$}, {\it J. Algebra }{\bf 319}
(2008), 1285--1304.

\re{Y} D. Yau, Enveloping algebras of hom-Lie algebras, {\it J. Gen.
Lie Theory Appl.} {\bf2} (2008), 95--108.

\re{S} Y. Su, Poisson brackets and structure of nongraded
Hamiltonian Lie algebras related to locally-finite derivations, {\it
Canad. J. Math.} {\bf55} (2003), 856--896.

%
%

\re{SX1} Y. Su, X. Xu, {Structure of divergence-free Lie algebras,}
{\it J. Algebra} {\bf 243} (2001), 557--595.

\re{SX2} Y. Su, X. Xu, {Central simple Poisson algebras}, {\it
Science in China A} {\bf 47} (2004), 245--263.

\re{SX3} Y. Su, X. Xu, Structure of contact Lie algebras related to
locally-finite derivations, {\it Manuscripta Math.} {\bf112} (2003),
231--257.

\re{SXZ} Y. Su, X. Xu, H. Zhang, {Derivation-simple Algebras and the
Structures of Lie Algebras of Witt type}, {\it J. Algebra} {\bf 233}
(2000), 23--58.

\re{SZ} Y. Su, K. Zhao, {Second cohomology group of generalized Witt
type Lie algebras and certain representations}, {\it Comm. Algebra }
{\bf 30} (2002), 3285--3309.

\re{SZh}    Y. Su, J. Zhou, Some representations of nongraded Lie
algebras of generalized Witt type, {\it J. Algebra} {\bf 246}
(2001), 721--738.

 \re{X} X. Xu,{ New generalized simple Lie algebras of
Cartan type over a field with characteristic 0}, {\it J. Algebra}
{\bf 224} (2000), 23--58.

%

\end{document}